\documentclass{article}
\usepackage[centertags]{amsmath}
\usepackage{amsfonts}
\usepackage{amssymb}
\usepackage{amsthm}
\usepackage{newlfont}
\usepackage{graphicx}
\newcommand \qq {\qquad}
\newcommand \q  {\quad}

\theoremstyle{plain}
\newtheorem{thm}{Theorem}[section]

\newtheorem{cor}[thm]{Corollary}

\theoremstyle{definition}

\numberwithin{equation}{section}

\begin{document}

\centerline{\textbf{ \begin{large} Constant Scalar Curvature of
Three Dimensional Surfaces
\end{large}}}
\centerline{\textbf{ \begin{large}  Obtained by the Equiform
Motion of a Sphere
\end{large}}}

\bigskip
\centerline{\begin{large} Fathi M. Hamdoon$^a$ , Ahmad T. Ali $^{b,}$ \footnote{Corresponding author.\\
E-mail address: atali71@yahoo.com (A. T. Ali).} and Rafael
L\'{o}pez$^c$
\end{large} }

\centerline{\footnotesize $^a$ Mathematics Department, Faculty of
Science, Al-Azhar University, Assiut, Egypt }
\centerline{\footnotesize $^b$ Mathematics Department, Faculty of
Science, Al-Azhar University, Cairo, Egypt }
\centerline{\footnotesize $^c$ Departamento de Geometr\'{\i}a y
Topolog\'{\i}a, Universidad de Granada, 18071 Granada, Spain } $ $
\hrule

\begin{flushleft}\textbf{Abstract}.\end{flushleft}

{\it In this paper we consider the equiform motion of a sphere in
Euclidean space $\mathbf{E}^7$. We study and analyze the
corresponding kinematic three dimensional surface under the
hypothesis that its scalar curvature $\mathbf{K}$ is constant. Under
this assumption, we prove that $|\mathbf{K}|<2$.}

\begin{flushleft}
\end{flushleft}

\begin{flushleft}
\end{flushleft}
\smallskip
\smallskip
\hrule

\section{Introduction}

An equiform transformation in the $n$-dimensional Euclidean space
$\mathbf{E}^n$ is an affine transformation whose linear part is
composed by an orthogonal transformation and a homothetical
transformation. Such an equiform transformation maps points
$\mathbf{x}\in\mathbf{E}^n$ according to the rule
\begin{equation}\label{equi1}
\mathbf{x}\longmapsto s {\cal A}\mathbf{x}+
\mathbf{d},\hspace*{.5cm} {\cal A}\in SO(n), s\in\mathbf{R}^+,
\mathbf{d}\in\mathbf{E}^n.
\end{equation}
The number $s$ is called the scaling factor. An equiform motion is
defined if the parameters of (\ref{equi1}), including $s$, are given
as functions of a time parameter $t$. Then a smooth one-parameter
equiform motion moves a point $\mathbf{x}$ via
$\mathbf{x}(t)=s(t){\cal A}(t) \mathbf{x}(t)+\mathbf{d}(t)$. The
kinematic corresponding to this transformation group is called
equiform kinematic. See  \cite{br,fhk}.

Under the assumption of the constancy of the scalar curvature,
kinematic surfaces obtained by the motion of a circle have been obtained in \cite{ah}. In
a similar context, one can consider hypersurfaces in  space forms
generated by  one-parameter family of spheres and having constant
curvature: \cite{cu,ja,lo,pa}.

In this paper we consider the equiform motions of a sphere
$\mathbf{k}_0$ in $\mathbf{E}^n$. The point paths of the sphere
generate a  $3$-dimensional surface, containing  the positions of
the starting sphere $\mathbf{k}_0$. The first order properties of
these surfaces for the points of these spheres have been studied for
arbitrary dimensions $n\geq 3$. We restrict our considerations to
dimension $n=7$ because, at any moment, the infinitesimal
transformations of the motion map the points of the sphere
$\mathbf{k}_0$ to the velocity vectors, whose end points will form
an affine image of $\mathbf{k}_0$ (in general a sphere) that span a
subspace $\mathbf{W}$ of $\mathbf{E}^n$ with $n\leq7$.

Let $\mathbf{x}(\theta,\phi)$ be a parameterization of
$\mathbf{k}_0$ and let $\mathbf{X}(t,\theta,\phi)$ be the resultant
3-surface by the equiform motion. We consider a certain position of
the moving space given by $t=0$, and we would like to obtain
information about the motion at  least during a certain period
around $t=0$  if we know its characteristics for one instant. Then
we  restrict our study to the properties of the motion for the limit
case $t\rightarrow 0$. A first choice is then approximate
$\mathbf{X}(t,\theta,\phi)$ by the first derivative of the
trajectories. Solliman, et al.  studied $3$-dimensional surfaces in
$\mathbf{E}^7$ generated by equiform motions of a sphere proving
that, in general, they are contained in a canal hypersurface
\cite{skhs}.

The purpose of this paper is to describe the kinematic surfaces
obtained by the motion of a sphere and whose scalar curvature
$\mathbf{K}$ is constant. As a consequence of our results, we prove:

\begin{quote} \emph{ A kinematic three-dimensional surface obtained by
the equiform motion of a sphere and with constant scalar curvature
$\mathbf{K}$ satisfies $|\mathbf{K}|<2$.}
\end{quote}
Moreover, we show the description of  the motion of such 3-surface
giving the equations that determine the kinematic geometry.

\section{The representation of a kinematic surface}

In two copies $\sum^{0}$, $\sum$ of Euclidean 7-space
$\mathbf{E}^7$,  we consider a unit sphere  $\mathbf{k}_0$ centered
at the origin of the 3-space $\varepsilon_0=[x_1x_2x_3]$ and
represented by
\begin{equation}\label{eq21}\mathbf{x}(\theta,\phi)=\Big(\cos(\theta)\cos(\phi),
\sin(\theta)\cos(\phi),\sin(\phi),0,0,0,0\Big)^{\text{T}},\,\,\,
\theta\in[0,2\pi],\,\,\,\phi\in[0,\pi].
\end{equation}
Under a one-parameter equiform motion of moving space $\sum^{0}$
with respect to a fixed space $\sum$ the general representation of
the motion of this surface  in $\mathbf{E}^7$ is given by
$$\mathbf{X}(t,\theta,\phi)=s(t)\mathbf{A}(t)\mathbf{x}(\theta,\phi)+\mathbf{d}(t),\qq t\in I\subset
\mathbf{R}.$$
 Here
$\mathbf{d}(t)=\Big(b_i(t)\Big)^{\text{T}}:i=1,2,...,7$ describes
the position of the origin of $\sum^{0}$ at time $t$,
$\mathbf{A}(t)=\Big(a_{ij}(t)\Big)^{\text{T}}:i,j=1,2,...,7$ is an
orthogonal matrix and $s(t)$ provides the scaling factor of the
moving system. With $s=\text{const.}\neq 0$ (sufficient to set
$s=1$), we have an ordinary Euclidean rigid body motion. For varying
$t$ and fixed $\mathbf{x}(\theta,\phi)$, equation (\ref{eq21}) gives
a parametric representation of the surface (or trajectory) of
$\mathbf{x}(\theta,\phi)$. Moreover, we assume that all involved
functions are at least of class $\mathbf{C}^1$. Using Taylor's
expansion up to the first order,  the representation of the motion
is given by
$$\mathbf{X}(t,\theta,\phi)=\Big[s(0)\mathbf{A}(0)+t\Big(\dot{s}(0)\mathbf{A}(0)+s(0)\dot{\mathbf{A}}(0)
\Big)\Big]\mathbf{x}(\theta,\phi)+\mathbf{d}(0)+t\dot{\mathbf{d}}(0),
$$ where $(.)$ denotes differentiation with respect to the time $t$.
Assuming that the moving frames $\sum^{0}$ and $\sum$ coincide at
the zero position ($t=0$),  we have
$$\mathbf{A}(0)=\mathbf{I},\q\, s(0)=1\q\, \text{and}\q\, \mathbf{d}(0)=0.
$$
Thus we have
$$\mathbf{X}(t,\theta,\phi)=\Big[\mathbf{I}+t\Big(s^{\prime}\mathbf{I}+\Omega
\Big)\Big]\mathbf{x}(\theta,\phi)+t\mathbf{d}^{\prime},$$
 where $\Omega=\dot{\mathbf{A}}(0)=(\omega_i),\,i=1,2,...,21$ is a skew symmetric matrix,
$s^{\prime}=\dot{s}(0)$, $\mathbf{d}^{\prime}=\dot{\mathbf{d}}(0)$
and all values of $s, b_i$ and their derivatives are computed at
$t=0$. With respect to these frames, the representation of the
motion up to the first order is
\begin{eqnarray}
\begin{pmatrix}
  \mathbf{X}_1 \\
   \mathbf{X}_2 \\
  \mathbf{X}_3 \\
   \mathbf{X}_4 \\
   \mathbf{X}_5 \\
   \mathbf{X}_6 \\
   \mathbf{X}_7
\end{pmatrix}=\begin{pmatrix}
 1+s^{\prime}t & \omega_1\,t & \omega_2\,t & \omega_3\,t & \omega_4\,t & \omega_5\,t & \omega_6\,t \\
  -\omega_1\,t & 1+s^{\prime}t & \omega_7\,t & \omega_8\,t & \omega_9\,t & \omega_{10}\,t & \omega_{11}\,t \\
  -\omega_2\,t & -\omega_7\,t & 1+s^{\prime}t & \omega_{12}\,t & \omega_{13}\,t & \omega_{14}\,t & \omega_{15}\,t \\
  -\omega_3\,t & -\omega_8\,t & -\omega_{12}\,t & 1+s^{\prime}t & \omega_{16}\,t & \omega_{17}\,t & \omega_{18}\,t \\
  -\omega_4\,t & -\omega_9\,t & -\omega_{13}\,t & -\omega_{16}\,t & 1+s^{\prime}t & \omega_{19}\,t & \omega_{20}\,t \\
  -\omega_5\,t & \omega_{10}\,t & -\omega_{14}\,t & -\omega_{17}\,t & -\omega_{19}\,t & 1+s^{\prime}t & \omega_{21}\,t \\
  -\omega_6\,t & -\omega_{11}\,t & -\omega_{15}\,t & -\omega_{18}\,t & -\omega_{20}\,t & -\omega_{21}\,t & 1+s^{\prime}t
\end{pmatrix}\,\begin{pmatrix}
 \cos(\theta)\cos(\phi) \\
 \sin(\theta)\cos(\phi) \\
 \sin(\phi) \\
  0 \\
  0 \\
  0 \\
  0
\end{pmatrix}+t\begin{pmatrix}
  b^{\prime}_{1} \\
  b^{\prime}_{2} \\
  b^{\prime}_{3} \\
  b^{\prime}_{4} \\
  b^{\prime}_{5} \\
  b^{\prime}_{6} \\
  b^{\prime}_{7}
\end{pmatrix},\nonumber
\end{eqnarray}
or in the equivalent form
\begin{eqnarray}\begin{pmatrix}
  \mathbf{X}_{1} \\
  \mathbf{X}_{2} \\
  \mathbf{X}_{3} \\
  \mathbf{X}_{4} \\
  \mathbf{X}_{5} \\
  \mathbf{X}_{6} \\
  \mathbf{X}_{7}
\end{pmatrix}&=&t\begin{pmatrix}
  b^{\prime}_{1} \\
  b^{\prime}_{2} \\
  b^{\prime}_{3} \\
  b^{\prime}_{4} \\
  b^{\prime}_{5} \\
  b^{\prime}_{6} \\
  b^{\prime}_{7}
\end{pmatrix}+\cos(\theta)\cos(\phi)\begin{pmatrix}
  1+s^{\prime}\,t \\
  -\omega_{1}\,t \\
  -\omega_{2}\,t \\
  -\omega_{3}\,t \\
  -\omega_{4}\,t \\
  -\omega_{5}\,t \\
  -\omega_{6}\,t
\end{pmatrix}+\sin(\theta)\cos(\phi)\begin{pmatrix}
  \omega_{1}\,t \\
  1+s^{\prime}\,t \\
  -\omega_{7}\,t \\
  -\omega_{8}\,t \\
  -\omega_{9}\,t \\
  -\omega_{10}\,t \\
  -\omega_{11}\,t
\end{pmatrix}+\sin(\phi)\begin{pmatrix}
  \omega_{2}\,t \\
   \omega_{7}\,t \\
  1+s^{\prime}\,t \\
  -\omega_{12}\,t \\
  -\omega_{13}\,t \\
  -\omega_{14}\,t \\
  -\omega_{15}\,t
\end{pmatrix}\nonumber\\
&=&t\,\vec{\mathbf{b}}+\cos(\theta)\cos(\phi)\,\vec{\mathbf{a}}_0
+\sin(\theta)\cos(\phi)\,\vec{\mathbf{a}}_1+\sin(\phi)\,\vec{\mathbf{a}}_2.\label{eq22}
\end{eqnarray}
For any fixed $t$ in in equation (\ref{eq22}), we generally get an
ellipsoid for $\theta\in[0,2\pi]$ and $\phi\in[0,\pi]$ centered at
the point
$t(b^{\prime}_1,b^{\prime}_2,b^{\prime}_3,b^{\prime}_4,b^{\prime}_5,b^{\prime}_6,b^{\prime}_7)$.
The latter ellipsoid turns to a 2-dimensional sphere if
$\vec{\mathbf{a}}_0$, $\vec{\mathbf{a}}_1$ and $\vec{\mathbf{a}}_2$
form an orthonormal basis. This gives the following conditions:
\begin{eqnarray}
\sum_{i=2}^{6}\omega_i\,\omega_{i+5}&=&\omega_1\,\omega_7-\sum_{i=3}^{6}\omega_i\,\omega_{i+9}=\omega_1\,\omega_2+
\sum_{i=8}^{11}\omega_i\,\omega_{i+4}=0,  \label{eq23}\\
\sum_{i=2}^{6}\omega_i^2&=&\sum_{i=7}^{11}\omega_i^2,\q\q\,\,
\omega_1^2+\sum_{i=3}^{6}\omega_i^2=\omega_7^2+\sum_{i=12}^{15}\omega_i^2.\label{eq24}
\end{eqnarray}

\section{ Scalar curvature of the kinematic surface }

In this section we shall compute the scalar curvature of 3-surfaces
in $\mathbf{E}^7$ generated by equiform motions of a sphere which
satisfies the conditions (\ref{eq23})-(\ref{eq24}). The tangents to
the parametric curves $t=\text{const}.$, $\theta=\text{const}.$ and
$\phi=\text{const}.$ at the zero position  are
$$\mathbf{X}_t=\Big[s^{\prime}\mathbf{I}+\Omega
\Big]\mathbf{x}+\mathbf{d}^{\prime},\,\,\,
\mathbf{X}_\theta=\Big[\mathbf{I}+\Big(s^{\prime}\mathbf{I}+\Omega\Big)t
\Big]\mathbf{x}_{\theta},\,\,\,
\mathbf{X}_\phi=\Big[\mathbf{I}+\Big(s^{\prime}\mathbf{I}+\Omega\Big)t
\Big]\mathbf{x}_{\phi},
$$
The first fundamental quantities of $\mathbf{X}(t,\theta,\phi)$ are
\begin{eqnarray*}
  \begin{cases}
g_{11}=\mathbf{X}_t\,\mathbf{X}_t^{\text{T}},\q\,
g_{12}=\mathbf{X}_\theta\,\mathbf{X}_t^{\text{T}},\q\,
g_{13}=\mathbf{X}_\phi\,\mathbf{X}_t^{\text{T}},\q\, & \text{ } \\
    g_{22}=\mathbf{X}_\theta\,\mathbf{X}_\theta^{\text{T}},\q\,
g_{23}=\mathbf{X}_\phi\,\mathbf{X}_\theta^{\text{T}},\q\,
g_{33}=\mathbf{X}_\phi\,\mathbf{X}_\phi^{\text{T}}. & \text{}
  \end{cases}
  \end{eqnarray*}
Under the conditions (\ref{eq23})-(\ref{eq24}), we obtain
\begin{eqnarray*}
g_{11}&=&\gamma+\alpha_{5}\cos(2\phi)+\alpha_{8}\sin(\phi)+2\cos(\phi)\Big[\cos(\phi)\Big(
\alpha_{4}\cos(2\theta)
+\alpha_{1}\sin(2\theta)\Big)\nonumber\\
&+&\sin(\theta)\Big(
\alpha_{7}+\alpha_{2}\sin(\phi)\Big)+\cos(\theta)\Big(
\alpha_{6}-2\alpha_{3}\sin(\phi)\Big)\Big],\\
g_{12}&=&\cos(\phi)\Big[
2\,t\,\cos(\phi)\Big(\alpha_{1}\cos(2\theta)-\alpha_{4}\sin(2\theta)\Big)-\omega_1\cos(\phi)\nonumber\\
&-&\sin(\theta)\Big[ t\Big(
\alpha_{6}-2\alpha_{3}\sin(\phi)\Big)+b^{\prime}_1+\omega_2\sin(\phi)\Big]\nonumber\\
&+&\cos(\theta)\Big[ t\Big(
\alpha_{7}+2\alpha_{2}\sin(\phi)\Big)+b^{\prime}_2+\omega_7\sin(\phi)\Big]\Big],\\
g_{13}&=&2t\cos(2\phi)\Big(
\alpha_{2}\sin(\theta)-\alpha_{3}\cos(\theta)\Big)-t\sin(2\phi)\Big(\alpha_{5}+\alpha_{4}\cos(2\theta)\nonumber\\
&+&\alpha_{1}\sin(2\theta)\Big)-\sin(\phi)\Big[
\Big(b^{\prime}_1+t\,\alpha_{6}\Big)\cos(\theta)+\Big(b^{\prime}_2+t\,\alpha_{7}\Big)\sin(\theta)\Big],\\
g_{22}&=&\cos^2(\phi)\Big[ 1+2t\,s^{\prime}+2t^2\Big(\delta-
\alpha_{4}\cos(2\theta)-\alpha_{6}\sin(2\theta)\Big)\Big],
\\
g_{23}&=&t^2\Big[ 2\cos^2(\phi)\Big(
\alpha_{2}\cos(\theta)+\alpha_{3}\sin(\theta)\Big)+\sin(2\phi)\Big(
\alpha_{4}\sin(2\theta)-\alpha_{1}\cos(2\theta)\Big)\Big],\\
g_{33}&=&1+2t\,s^{\prime}+t^2\Big[\gamma-\beta-\alpha_{5}\cos(2\phi)+2\sin^2(\phi)\Big(\alpha_{4}\cos(2\theta)+\alpha_{1}\sin(2\theta)\Big)\nonumber\\
&+&2\sin(2\phi)
\Big(\alpha_{3}\cos(\theta)-\alpha_{2}\sin(\theta)\Big)\Big],
\end{eqnarray*}
where
\begin{eqnarray*}
\begin{cases}
\alpha_{1}=\dfrac{1}{2}\Big[\sum_{i=2}^6
\omega_{i}\,\omega_{i+5}\Big], & \\
\alpha_{2}=\dfrac{1}{2}\Big[\omega_1\,\omega_2+\sum_{i=8}^{11}
\omega_{i}\,\omega_{i+4}\Big], & \\
\alpha_{3}=\dfrac{1}{2}\Big[\omega_1\,\omega_7-\sum_{i=3}^{6}
\omega_{i}\,\omega_{i+9}\Big],& \\ \alpha_{4}=\dfrac{1}{4}\Big[\sum_{i=2}^{6}\Big(\omega_{i}^2-\omega_{i+5}^2\Big)\Big], & \\
\alpha_{5}=
\dfrac{1}{4}\Big[\omega_1^2-2\,\omega_2^2-2\,\omega_7^2+\sum_{i=1}^{11}\omega_{i}^2-2\Big(\sum_{i=12}^{15}\omega_{i}^2\Big)\Big],
& \\
\alpha_{6}=b^{\prime}_1\,s^{\prime}-\sum_{i=2}^{7}
b^{\prime}_{i}\,\omega_{i-1}, & \\
\alpha_{7}=b^{\prime}_1\,\omega_1+b^{\prime}_2\,s^{\prime}-\sum_{i=3}^{7}
b^{\prime}_{i}\,\omega_{i+4}, & \\
\alpha_{8}=2\Big[b^{\prime}_1\,\omega_2+b^{\prime}_2\,\omega_7+
b^{\prime}_3\,s^{\prime}-\sum_{i=4}^{7}
b^{\prime}_{i}\,\omega_{i+8}\Big], & \\
\beta=\sum_{i=1}^{7} b_{i}^{\prime2}, & \\
\gamma=\beta+s^{\prime2}+\dfrac{1}{4}\Big[2(\omega_1^2+\omega_2^2+\omega_7^2)+\sum_{i=2}^{15}
\omega_{i}^2+\sum_{i=12}^{15} \omega_{i}^2\Big], & \\
\delta=\dfrac{1}{4}\Big[2\,(s^{\prime2}+\omega_1^2)+\sum_{i=2}^{11}
\omega_{i}^2\Big].
\end{cases}
\end{eqnarray*}
The conditions (\ref{eq23})-(\ref{eq24}) lead to the following
relations
$$\alpha_{1}=\alpha_{2}=\alpha_{3}=\alpha_{4}=\alpha_{5}=0,\q
\gamma=\beta+2\delta.
$$
In order to calculate the scalar curvature, we need to compute the
Christoffel symbols of the second kind, which are defined as
\begin{equation}\label{eq311}
\Gamma_{ij}^{l}=\dfrac{1}{2}g^{lm}\Big[\frac{\partial
g_{im}}{\partial x_{j}}+\frac{\partial g_{jm}}{\partial
x_{i}}-\frac{\partial g_{ij}}{\partial x_{m}}\Big],
\end{equation}
where $i,j,l$ are indices that take the values $1,2,3$, $x_1=t,
x_2=\theta, x_3=\phi$, and $\Big(g^{lm}\Big)$ is the inverse matrix
of $\Big(g_{ij}\Big)$. Then the scalar curvature of the surface
$\mathbf{X}(t,\theta,\phi)$ is
$$\mathbf{K}(t,\theta,\phi)=g^{ij}\Big[\frac{\partial
\Gamma_{ij}^l}{\partial x_{l}}-\frac{\partial
\Gamma_{il}^l}{\partial
x_{j}}+\Gamma_{ij}^l\,\Gamma_{lm}^m-\Gamma_{il}^m\,\Gamma_{jm}^l\Big].$$
At the zero position ($t=0$), the scalar curvature of
$\mathbf{X}(t,\theta,\phi)$ is given by
\begin{eqnarray}
\mathbf{K}=\mathbf{K}(0,\theta,\phi)&=&\dfrac{P\Big(\cos(n_1\theta\pm
m_1\phi),\sin(n_1\theta\pm m_1\phi)\Big)}{
Q\Big(\cos(n_2\theta\pm m_2\phi),\sin(n_2\theta\pm m_2\phi)\Big)}.\label{eq31}
\end{eqnarray}
This quotient writes then as
\begin{equation}\label{eq32}
P\Big(\cos(n_1\theta\pm m_1\phi),\sin(n_1\theta\pm
m_1\phi)\Big)-\mathbf{K}\,Q\Big(\cos(n_2\theta\pm
m_2\phi),\sin(n_2\theta\pm m_2\phi)\Big)=0.
\end{equation}
The assumption on the constancy of the scalar curvature $\mathbf{K}$
implies that  equation (\ref{eq32}) is a linear combination of the
functions $\cos(n\,\theta\pm m\,\phi)$, $\sin(n\,\theta\pm
m\,\phi)$. Because these functions are  independent linearly, the
corresponding coefficients must  vanish. Throughout this work, we
have employed the Mathematica programm in order to compute the
explicit expressions of these coefficients.

{\bf Assumption.} Without loss of generality, we assume that the two
conditions (\ref{eq23})-(\ref{eq24}) are satisfied and there are no
translation motions in the plane which contain the starting sphere,
i.e.,
$$
b^{\prime}_1=b^{\prime}_2=b^{\prime}_3=0. $$

\subsection{ Kinematic surfaces with zero scalar curvature}

We assume that $\mathbf{K}=0$. From the expression (\ref{eq31}), we
have
\begin{eqnarray*}
P\Big(\cos(n_1\theta&\pm &m_1\phi),\sin(n_1\theta\pm m_1\phi)\Big)\nonumber\\
&=&\sum_{i=0}^{12}\sum_{j=-12}^{12}\Big(A_{i,j}\cos(i\,\theta+j\,\phi)+B_{i,j}\sin(i\,\theta+j\,\phi)\Big)=0.
\end{eqnarray*}
In this case, a straightforward computation shows that the
coefficients of $\cos(12\phi)$, $\cos(6\theta+12\phi)$ and
$\sin(6\theta+12\phi)$ are
\begin{eqnarray*}
A_{0,12}&=&\frac{3}{8192}\Big[16\,\omega_1^6-120\,\omega_1^4\,\Big(\omega_2^2+\omega_7^2\Big)+9\,\omega_1^2\,
\Big(\omega_2^2+\omega_7^2\Big)^2-5\,\Big(\omega_2^2+\omega_7^2\Big)^3\Big] \nonumber\\
A_{6,12}&=&\frac{3}{32768}\Big[\omega_2^6-15\,\omega_2^4\,\omega_7^2+15\,\omega_2^2\,\omega_7^4-\omega_7^6\Big] \nonumber\\
B_{6,12}&=&\frac{3}{16384}\,\omega_2\,\omega_7\Big(3\,\omega_2-\omega_7^2\Big)\Big(\omega_2^2-3\,\omega_7^2\Big).
\nonumber
\end{eqnarray*}
By solving the three equations $A_{0,12}=0$, $A_{6,+12}=0$ and
$B_{6,+12}=0$, we get
$$ \omega_1=\omega_2=\omega_7=0.$$
Then
\begin{eqnarray*}
B_{0,9}&=&\frac{3}{256}\,\alpha_{8}\,\Big[\alpha_{8}^2-6\,\Big(\alpha_6^2+\alpha_{7}^2\Big)\Big]\nonumber\\
A_{3,9}&=&\frac{3}{256}\,\alpha_{6}\,\Big(3\,\alpha_{7}^2-\alpha_{6}^2\Big)\nonumber\\
B_{3,9}&=&\frac{3}{256}\,\alpha_{7}\,\Big(\alpha_{7}^2-3\,\alpha_{6}^2\Big).\nonumber
\end{eqnarray*}
The three equations $B_{0,9}$, $A_{3,9}=0$ and $B_{3,9}=0$ imply
$$\alpha_{6}=\alpha_{7}=\alpha_8=0. $$
 From these values,
equation $A_{0,6}=0$ leads to
$$
\Big(\beta+2\,\delta\Big)\Big(\beta+s^{\prime2}-2\,\delta\Big)=0.
$$
 It is worth to point out  that the quantities
$\beta$ and $\delta$ are positive and thus  we obtain the following
condition
$$\delta=\frac{1}{2}\Big(\beta+s^{\prime2}\Big).$$
At this time, the explicit computations of coefficients imply that
all $A_{i,j}$ and $B_{i,j}$ are equal zero. So, we have the
following

\begin{thm} A kinematic 3-surface in $\mathbf{E}^7$ foliated by spheres and with zero constant scalar curvature
satisfies
\begin{eqnarray*}
  \omega_1&=&\omega_2=\omega_7=0,   \\
    \sum_{i=4}^{7}b^{\prime}_i\,\omega_{i-1}&=&
\sum_{i=4}^{7}b^{\prime}_i\,\omega_{i+4}=
\sum_{i=4}^{7}b^{\prime}_i\,\omega_{i+8}=0,  \\
\sum_{i=3}^{6}\omega_i^2&=&\sum_{i=4}^{7}b^{\prime2}_i.
    \end{eqnarray*}
\end{thm}

\subsection{ Kinematic surfaces with non-zero constant scalar curvature}

We assume that the kinematic 3-surface has  constant scalar
curvature $\mathbf{K}\not= 0$. From (\ref{eq311}), we have
\begin{eqnarray*}
P\Big(\cos(n_1\theta&\pm &m_1\phi),\sin(n_1\theta\pm m_1\phi)\Big)
-\mathbf{K}\,Q\Big(\cos(n_2\theta\pm m_2\phi),\sin(n_2\theta\pm m_2\phi)\Big)\nonumber\\
&=&\sum_{j=-12}^{12}\sum_{i=0}^{12}\Big(A_{i,j}\cos(i\,\theta+j\,\phi)+B_{i,j}\sin(i\,\theta+j\,\phi)\Big)=0.
\end{eqnarray*}
In this case, a straightforward computation shows that the
coefficients of $\cos(12\phi)$, $\cos(12\theta+6\phi)$ and
$\sin(12\theta+6\phi)$ are
\begin{eqnarray*}
A_{0,12}&=&\frac{1}{16384}\,\Big(\mathbf{K}+6\Big)\Big[16\omega_1^6-120\,\omega_1^4\,\Big(\omega_2^2+\omega_7^2\Big)+
90\,\omega_1^2\,\Big( \omega_2^2+\omega_7^2\Big)^2-5\,\Big(
\omega_2^2+\omega_7^2\Big)^3\Big],\nonumber\\
A_{6,12}&=&\frac{1}{65536}\,\Big(\mathbf{K}+6\Big)\Big[\omega_2^6-15\,\omega_2^4\,\omega_7^2+15\omega_2^2\,\omega_7^4
-\omega_7^6\Big],\nonumber\\
B_{6,12}&=&\frac{1}{32768}\,\omega_2\,\omega_7\,\Big(\mathbf{K}+6\Big)\Big(3\,\omega_2-\omega_2^2\Big)\Big(\omega_2^2-3\,
\omega_7^2\Big). \nonumber
\end{eqnarray*}
We consider the three equations  $A_{0,12}=0$, $A_{6,12}=0$ and
$B_{6,12}=0$. From here, we discuss  two possibilities:
$\mathbf{K}=-6$ and $\omega_1=\omega_2=\omega_7=0$.

\begin{enumerate}
\item {\bf{Case  $\mathbf{K}=-6$}}. A  computation of coefficients
yields
\begin{eqnarray*}
A_{5,11}&=&\frac{1}{2048}\Big[\alpha_{6}\Big(\omega_2^4-6\,\omega_2^2\,\omega_7^2+\omega_7^4\Big)
-4\alpha_{7}\,\omega_2\,\omega_7\Big(\omega_2^2-\omega_7^2\Big)\Big]=0 \nonumber\\
B_{5,11}&=&\frac{1}{2048}\Big[\alpha_{7}\Big(\omega_2^4-6\,\omega_2^2\,\omega_7^2+\omega_7^4\Big)
+4\alpha_{6}\,\omega_2\,\omega_7\Big(\omega_2^2-\omega_7^2\Big)\Big]=0.
\nonumber
\end{eqnarray*}
We consider two cases: $\alpha_{6}=\alpha_{7}=0$ and
$\omega_2=\omega_7=0$.

{\bf{Case (1):}} We assume $\alpha_{6}=\alpha_{7}=0$. The
computation of coefficients leads to
\begin{eqnarray*}
B_{0,11}&=&\frac{1}{2048}\,\alpha_8\,\Big[8\,\omega_1^4-24\,\omega_1^2\,\Big(\omega_2^2+\omega_7^2\Big)
+3\,\Big(\omega_2^2+\omega_7^2\Big)^2\Big]=0, \nonumber\\
A_{4,11}&=&\frac{1}{512}\,\alpha_{8}\,\omega_2\,\omega_7\Big(\omega_7^2-\omega_2^2\Big)=0, \nonumber\\
B_{4,11}&=&\frac{1}{2048}\,\alpha_{8}\,\Big(\omega_2^4-6\,\omega_2^2\,\omega_7^2+\omega_7^4\Big)=0,
\nonumber
\end{eqnarray*}
which  implies two subcases: $\alpha_{8}=0$ and
$\omega_1=\omega_2=\omega_7=0$.

{\bf{Subcase (1.1):}} If $\alpha_{8}=0$, then we have
\begin{eqnarray*}
A_{0,10}&=&\frac{1}{1024}\,\Big[8\,\omega_1^2-\omega_1^2\Big(\omega_2^2+\omega_7^2\Big)+3\,\Big(\omega_2^2+\omega_7^2\Big)^2
\Big]
\Big(\beta+s^{\prime2}+6\,\delta-\omega_1^2-\omega_2^2-\omega_7^2\Big) \nonumber\\
A_{4,10}&=&\frac{1}{2048}\,\omega_2\,\omega_7\,\Big(\omega_2^4-6\,\omega_2^2\,\omega_7^2+\omega_7^4\Big)
\Big(\beta+s^{\prime2}+6\,\delta-\omega_1^2-\omega_2^2-\omega_7^2\Big) \nonumber\\
B_{4,10}&=&\frac{1}{512}\,\omega_2\,\omega_7\,\Big(\omega_7^2-\omega_2^2\Big)
\Big(\beta+s^{\prime2}+6\,\delta-\omega_1^2-\omega_2^2-\omega_7^2\Big).
\nonumber
\end{eqnarray*}
The last term in the above three equations is not  zero because
$$
\beta+s^{\prime2}+6\,\delta-\omega_1^2-\omega_2^2-\omega_7^2=\sum_{i=4}^7
b^{\prime2}_i+\omega_2^2+\sum_{i=8}^{11}\omega_i^2+2\Big[2\,s^{\prime2}+\omega_1^2+\sum_{i=3}^6\omega_i^2\Big]>0.
$$
The three equations $A_{0,10}=0$, $A_{4,10}=0$ and $B_{4,10}=0$ lead
$\omega_1=\omega_{2}=\omega_{7}=0$. Now, the coefficient $A_{0,6}$
must equal zero, that is,
$$\Big(\beta+2\,\delta\Big)^2\,\Big(2\,\beta+8\,\delta-s^{\prime2}\Big)=0,
$$ contradiction.

{\bf{Subcase (1.2):}} If $\omega_1=\omega_2=\omega_7=0$ and
$\alpha_{8}\neq0$, the equation $B_{0,9}=0$ implies that
$\alpha_{8}=0$: contradiction.

{\bf{Case (2):}} If $\omega_2=\omega_7=0$ and
$\alpha_{6},\alpha_7\neq0$, the computation of coefficients yields
\begin{eqnarray*}
A_{5,11}=\frac{1}{2048}\,\alpha_{6}\,\omega_1^4=0,\nonumber\\
B_{5,11}=\frac{1}{2048}\,\alpha_{7}\,\omega_1^4=0. \nonumber
\end{eqnarray*}
Because $\alpha_{6}\neq0$ and $\alpha_{7}\neq0$, we conclude
$\omega_1=0$. New computations give
\begin{eqnarray*}
A_{3,9}=\frac{9}{256}\alpha_{6}\,\Big(\alpha_{6}^2-3\,\alpha_{7}^2\Big),\nonumber\\
B_{3,9}=\frac{9}{256}\alpha_{7}\,\Big(3\,\alpha_{6}^2-\alpha_{7}^2\Big),\nonumber
\end{eqnarray*}
By solving the equations $A_{3,9}=0$ and $B_{3,9}=0$, we get
$\alpha_{6}=\alpha_{7}=0$: contradiction.

\begin{cor} There are not kinematic 3-surfaces in
$\mathbf{E}^7$ foliated by spheres and with scalar curvature
$\mathbf{K}$ equal $-6$.
\end{cor}

\item {\bf{Case $
\omega_1=\omega_2=\omega_7=0$ and $\mathbf{K}\neq-6$.}}

A computation of the coefficients yields
\begin{eqnarray*}
B_{0,9}&=&\frac{1}{256}\,\alpha_{8}\,\Big[\alpha_{8}^2-6\,\Big(\alpha_{6}^2+\alpha_7^2\Big)\Big]\Big(2\,\mathbf{K}+3\Big)=0,\nonumber\\
A_{3,9}&=&\frac{1}{256}\,\alpha_{6}\,\Big(3\,\alpha_{7}^2-\alpha_{6}^2\Big)\Big(2\,\mathbf{K}+3\Big)=0,\nonumber\\
B_{3,9}&=&\frac{1}{256}\,\alpha_{7}\,\Big(\alpha_{7}^2-3\,\alpha_{6}^2\Big)\Big(2\,\mathbf{K}+3\Big)=0,\nonumber
\end{eqnarray*}
which gives two cases: $\mathbf{K}=-\dfrac{3}{2}$ or
$\alpha_{7}=\alpha_{6}=\alpha_8=0$.

{\bf{Case (1):}} Assume $\mathbf{K}=-\dfrac{3}{2}$. Now, we obtain
\begin{eqnarray*}
A_{0,8}&=&\frac{1}{64}\,\Big[\alpha_{8}^2-\Big(\alpha_{6}^2-\alpha_7^2\Big)\Big]\Big(6\,\delta-\beta-2\,s^{\prime2}\Big),\nonumber\\
A_{2,8}&=&\frac{1}{64}\,\Big(\alpha_{6}^2-\alpha_{7}^2\Big)\Big(6\,\delta-\beta-2\,s^{\prime2}\Big),\nonumber\\
B_{2,8}&=&\frac{1}{32}\,\alpha_{7}\alpha_{6}\Big(6\,\delta-\beta-2\,s^{\prime2}\Big).\nonumber
\end{eqnarray*}
Solving the three equations $A_{0,8}=0$, $A_{2,8}=0$ and
$B_{2,8}=0$, we find two cases:
$\delta=\dfrac{\beta+2\,s^{\prime2}}{6}$ and
$\alpha_{6}=\alpha_{7}=\alpha_8=0$.

{\bf{Case (1.1):}} If $\delta=\dfrac{\beta+2\,s^{\prime2}}{6}$, we
obtain
$$A_{0,4}=\frac{1}{72}\,(2\beta+s^{\prime2})\,\Big[
4(2\,\beta+s^{\prime2})^2-9\,\Big(\alpha_{8}^2+4(\alpha_{6}^2+\alpha_{7}^2)\Big)\Big],$$
which leads the following condition
$$4(2\,\beta+s^{\prime2})^2=9\,\Big(\alpha_{8}^2
+4(\alpha_{6}^2+\alpha_{7}^2)\Big). $$
At this point, all
coefficients $A_{i,j}$ and $B_{i,j}$ are equal zero.

{\bf{Case (1.2):}} Assume $\alpha_{6}=\alpha_{7}=\alpha_8=0$  and
$\delta\neq\dfrac{\beta+2\,s^{\prime2}}{6}$. The coefficient
$A_{0,6}$ is
$A_{0,6}=\dfrac{1}{32}\,\Big(\beta+2\,\delta\Big)\Big[14\,\delta-\beta-4\,s^{\prime2}\Big]
$. From $A_{0,6}=0$, we conclude
$$\delta=\dfrac{\beta+4\,s^{\prime2}}{14}.$$
As a consequence, all coefficients $A_{i,j}$ and $B_{i,j}$,
$i=1,2,...,12$, $j=-12,...,12$ are zero.

From the above reasonings, it follows the next

\begin{thm}\label{th-33} A kinematic 3-surface in $\mathbf{E}^7$ foliated by spheres and
with $\mathbf{K}=-\dfrac{3}{2}$ satisfies
$\omega_1=\omega_2=\omega_7=0$ and one of the next pairs of
equations:
\begin{eqnarray*}
s^{\prime2}+3\sum_{i=3}^6\omega_i^2&=&\sum_{i=4}^7 b^{\prime2}_i,
\\
4\Big[s^{\prime2}+2\sum_{i=4}^7 b^{\prime2}_i\Big]&=&9\Big[ \Big(
\sum_{i=4}^{7}b^{\prime}_i\omega_{i+8}\Big)^2+4\Big(
\sum_{i=4}^{7}b^{\prime}_i\omega_{i-1}\Big)^2+4\Big(\sum_{i=4}^{7}
b^{\prime}_i\omega_{i+4}\Big)^2\Big],
\end{eqnarray*}
or
\begin{eqnarray*}
\sum_{i=4}^{7}b^{\prime}_i\,\omega_{i-1}&=&
\sum_{i=4}^{7}b^{\prime}_i\,\omega_{i+4}=
\sum_{i=4}^{7}b^{\prime}_i\,\omega_{i+8}=0,\\
3\,s^{\prime2}&+&7\,\sum_{i=3}^6\,\omega_i^2=\sum_{i=4}^7\,b^{\prime2}_i.
\end{eqnarray*}
\end{thm}

{\bf{Case (2):}} Assume $\alpha_{6}=\alpha_{7}=\alpha_8=0$ and
$\mathbf{K}\neq-\dfrac{3}{2}$.  In this case we obtain
$$A_{0,6}=\frac{1}{16}\Big(\beta+2\,\delta\Big)^2\Big[\mathbf{K}\Big(\beta+2\,\delta\Big)
-2\Big(2\,\delta-\beta-s^{\prime2}\Big)\Big]=0, $$
 which yields
\begin{equation}\label{eq33}
\mathbf{K}=\frac{2(2\,\delta-\beta-s^{\prime2})}{\beta+2\delta}.
\end{equation}
From here, all coefficients $A_{i,j}$ and $B_{i,j}$ are equal zero.
So, we have the following

\begin{thm}\label{th-34} A kinematic 3-surface in $\mathbf{E}^7$ foliated by spheres and
with constant scalar curvature
$$\mathbf{K}=\frac{2\Big[\sum_{i=3}^6\,\omega_i^2-\sum_{i=4}^7\,b^{\prime2}\Big]}{
s^{\prime2}+\sum_{i=3}^6\,\omega_i^2+\sum_{i=4}^7\,b^{\prime2}}$$
satisfies
\begin{eqnarray*}
\omega_1&=&\omega_2=\omega_7=0,\\
\sum_{i=4}^{7}b^{\prime}_i\,\omega_{i-1}&=&
\sum_{i=4}^{7}b^{\prime}_i\,\omega_{i+4}=
\sum_{i=4}^{7}b^{\prime}_i\,\omega_{i+8}=0.
\end{eqnarray*}
\end{thm}

From the expression (\ref{eq33}), we can write the quantity $\beta$
in the form
$$\beta=\frac{2\Big[\Big(2-\mathbf{K}\Big)\delta-s^{\prime2}\Big]}{\mathbf{K}+2}.$$
As $\beta$ is positive,  we have two cases:
\begin{enumerate}
\item Case $\mathbf{K}+2<0$ and
$\Big(2-\mathbf{K}\Big)\delta-s^{\prime2}<0$. This implies
$$\mathbf{K}<-2,\q \text{and} \q
\mathbf{K}>\dfrac{2\,\delta-s^{\prime2}}{\delta}=\dfrac{
2\sum_{i=3}^6\,\omega_i^2}{s^{\prime2}+\sum_{i=3}^6\,\omega_i^2}>0,$$
which is a contradiction.

\item Case $\mathbf{K}+2>0$ and $\Big(2-\mathbf{K}\Big)\delta-s^{\prime2}>0$. This gives
 the following condition for $\mathbf{K}$:
$$ -2<\mathbf{K}<\dfrac{2\delta-s^{\prime2}}{\delta}=\dfrac{
2\sum_{i=3}^6\omega_i^2}{s^{\prime2}+\sum_{i=3}^6\,\omega_i^2}<2.$$
\end{enumerate}
\end{enumerate}

As consequence of Theorems \ref{th-33} and \ref{th-34}, we have the
next statement, which was established in the Introduction:

\begin{cor}  A kinematic three-dimensional surface in $\mathbf{E}^7$ obtained by
the equiform motion of a sphere and with constant scalar curvature
$\mathbf{K}$ satisfies $|\mathbf{K}|<2$.
\end{cor}



\begin{thebibliography}{99}

\bibitem{ah}  Abdel-All, N. H. and  Hamdoon, F. M.; \emph{Cyclic surfaces
in $E^5$ generated by equiform motions}, J. Geom. {\bf 79} (2004)
1-11.

\bibitem{br} Bottema O. and Roth B.; Theoretical kinematic, Dover publications
Inc., New York, 1990.

\bibitem{cu} Castro I. and Urbano F.; {\it On a minimal Lagrangian
submanifold of $\mathbf{C}^n$ foliated by spheres}, Michigan Math.
J., {\bf{46}}, (1999), 71-82.

\bibitem{fhk} Farin G., Hoschek J. and Kim M., The handbook of computer aided
geometric design, North-Holland, Amsterdam, 2002.

\bibitem{ja} Jagy W.; {\it Sphere foliated constant mean curvature submanifolds},
Rocky Mount. J. Math., {\bf{28}}, (1998), 983-1015.

\bibitem{lo}  L\'{o}pez R.; {\it Cyclic hypersurfaces of constant curvature },
Adv. Stud. Pure Math., {\bf{34}}, (2002), 185-199.

\bibitem{pa} Park S.H.; {\it Sphere foliated minimal and constant
mean curvature hypersurfaces in space forms Lorentz-Minkowski
space}, Rocky Mount. J. Math., {\bf{32}}, (2002), 1019-1044.

\bibitem{skhs} Solliman M., Khater A. Hamdoon F. and Solouma E.; {\it Three
dimensional surfaces foliated by two dimensional spheres }, J.
Egyptian Math. Soc., {\bf{15}}, (2007), 101-110.

\end{thebibliography}
\end{document}